\documentclass[11pt,dutch]{article}

\usepackage{graphicx}
\usepackage{tocbibind}
\usepackage{graphicx}
\usepackage{amsfonts}
\usepackage{amssymb}
\usepackage{amsthm}
\usepackage{amsmath}
\usepackage{yfonts}
\usepackage{eufrak}
\usepackage[all]{xy}
\usepackage{makeidx}
\makeindex

\theoremstyle{plain}

\newtheorem{theorem}{Theorem}[section]
\newtheorem{lemma}[theorem]{Lemma}
\newtheorem{proposition}[theorem]{Proposition}
\newtheorem{corollary}[theorem]{Corollary}
\newtheorem*{property-P}{Property P}

\theoremstyle{definition}

\newtheorem{definition}[theorem]{Definition}

\newtheorem{remark}[theorem]{Remark}

\newtheorem{example}[theorem]{Example}

\renewcommand{\to}{\longrightarrow}

\newcommand{\del}{\partial}
\newcommand{\im}{\ensuremath{\mathsf{Im\,}}}
\renewcommand{\im}{\ensuremath{\mathrm{Im}}}
\newcommand{\coker}{\ensuremath{\mathsf{Coker\,}}}
\newcommand{\Cok}{\ensuremath{\mathrm{Cok}}}
\renewcommand{\ker}{\ensuremath{\mathsf{Ker\,}}}

\newcommand{\Bc}{\ensuremath{\mathcal{B}}}
\newcommand{\Ac}{\ensuremath{\mathcal{A}}}
\newcommand{\Pc}{\ensuremath{\mathcal{P}}}

\newcommand{\Ab}{\ensuremath{\mathsf{Ab}}}
\newcommand{\ab}{\ensuremath{\mathrm{ab}}}

\newcommand{\Xc}{\ensuremath{\mathcal{X}}}

\newcommand{\Z}{\ensuremath{\mathbb{Z}}}
\newcommand{\N}{\ensuremath{\mathbb{N}}}

\newcommand{\Gp}{\ensuremath{\mathsf{Gp}}}

\newcommand{\Ext}{\ensuremath{\mathsf{Ext}}}
\newcommand{\Extn}{\ensuremath{\mathsf{Ext}^{n}\!}}

\newcommand{\Reg}{\ensuremath{\mathsf{Reg}}}

\newcommand{\Arr}{\ensuremath{\mathsf{Arr}}}

\newcommand{\Ob}{\ensuremath{\mathsf{Ob}}}
\newcommand{\Arrn}{\ensuremath{\mathsf{Arr}^{n}\!}}

\newcommand{\CExt}{\ensuremath{\mathsf{CExt}}}
\newcommand{\TExt}{\ensuremath{\mathsf{TExt}}}

\newcommand{\Ec}{\ensuremath{\mathcal{E}}}

\newcommand{\op}{\ensuremath{\mathrm{op}}}

\newbox\skewpullbackbox
\setbox\skewpullbackbox=\hbox{\xy 0;<1mm,0mm>: \POS(4,0)\ar@{-} (-4,0) \ar@{-} (8,4)
\endxy}

\newbox\ksewpullbackbox
\setbox\ksewpullbackbox=\hbox{\xy 0;<1mm,0mm>: \POS(0,-8)\ar@{-} (0,-4) \ar@{-} (4,-4)
\endxy}

\newbox\pullbackbox
\setbox\pullbackbox=\hbox{\xy 0;<1mm,0mm>: \POS(4,0)\ar@{-} (0,0) \ar@{-} (4,4)
\endxy}
\def\pullback{\copy\pullbackbox}

\newbox\pushoutbox
\setbox\pushoutbox=\hbox{\xy 0;<1mm,0mm>: \POS(0,4)\ar@{-} (0,0) \ar@{-} (4,4)
\endxy}

\begin{document}

\title{Higher central extensions and Hopf formulae}

\author{Tomas Everaert\footnote{The author is a postdoctoral fellow with FWO-Vlaanderen.}}
\date{}

\maketitle

\setcounter{section}{-1}
\begin{abstract}
\noindent
Higher extensions and higher central extensions, which are of importance to non-abelian homological algebra, are studied, and some fundamental properties are proven.
As an application, a direct proof of the invariance of the higher Hopf formulae is obtained.
\end{abstract}

\section{Introduction}
When Hopf discovered in the early $1940$'s a purely algebraic method to construct the second homology group of a pathwise connected aspherical complex from the fundamental group of the complex \cite{Hopf}, it marked the birth of what is now called homological algebra. Hopf proved that, when $A$ is the fundamental group of the considered complex, and $p\colon P\to A$ is a projective presentation of $A$ (i.e.\ a surjective homomorphism onto $A$ with a free domain), then the second (integral) homology group of the complex is isomorphic to the quotient
\begin{equation}\label{Hopfformule}
\frac{[P,P] \cap K[p]}{[K[p], P]},
\end{equation}
where $K[p]$ denotes the kernel of $p$ and $[\cdot,\cdot]$ the group commutator. In particular, he proved that the construction is independent of the particular choice of projective presentation of $A$. Thanks to Hopf's result it became meaningful to speak of the homology group of a \emph{group} rather than of a \emph{complex} (or, more generally, a \emph{topological space}). 

Constructions for the higher homology groups were found soon after: Eilenberg and MacLane, Hopf, Freudenthal and Eckmann all more or less independently came up with solutions that resulted in the now classic definition via chain resolutions. (See \cite{MacLane:Origins} for a detailed historical account.) Later, another now classic construction, via simplicial resolutions, was discovered (see, e.g.\ \cite{Barr-Beck}). It is interesting to note that, contrary to what one might expect, neither of these two coincides with Hopf's original construction, in the particular case of the second homology group. In fact, it was not until $1989$ that similar constructions were found for the higher homology groups: in their article \cite{Brown-Ellis}, Brown and Ellis obtained a ``Hopf formula'' for the $n$-th integral homology group, for any $n\geq 3$ (see also \cite{Ellis:Hopf}). 

To give an idea, let us have a look at Brown and Ellis's Hopf formula for the third homology of a group. Consider a group $A$. Instead of representing $A$ by a single surjective homomorphism, represent it by a commutative square of surjective homomorphisms 
\begin{equation}\label{doublepresentation}
\vcenter{\xymatrix{
F_{2} \ar[r]^{f_2} \ar[d]_{f'_2} & F_{1} \ar[d]^{f_1} \\
F_{0} \ar[r]_{f_0} & A}}
\end{equation}
 with the following properties: the unique factorization $a\colon F_2\to P$ to the pullback $P=F_{0}\times_A F_1$ is surjective, and $F_0$, $F_1$ and $F_{2}$ are free groups. Brown and Ellis discovered that the third (integral) homology group of $A$ is isomorphic to the following quotient, where the denominator is an internal product of subgroups of $A$:
\begin{equation}\label{Hopf3}
\frac{[F_2,F_2]\cap K[f_2]\cap K[f'_2]}{[K[f_2],K[f'_2]] \cdot [K[f_2]\cap K[f'_2], F_2]}.
\end{equation}

Recall that a surjective homomorphism $f\colon A\to B$ is also called an \emph{extension} (of $B$). Commutative squares of extensions, such as \ref{doublepresentation}, in which also the factorization $a\colon F_2\to P$ to the pullback $P=F_{0}\times_A F_1$ is an extension (but without any freeness condition) are called \emph{double extensions}. They have been considered  by Janelidze in his study of \emph{double central extensions} \cite{Janelidze:Double}, as an application of Categorical Galois Theory \cite{Janelidze:Pure}. In his talk \cite{Janelidze:Hopf-talk} Janelidze continued his observations from \cite{Janelidze:Double} by considering \emph{$n$-fold extensions} and \emph{$n$-fold central extensions} for \emph{any} $n\geq 1$. Both are particular kinds of commutative $n$-dimensional cubes of extensions. Higher extensions are of interest: indeed, Brown and Ellis's construction of the $(n+1)$-st homology of a group $A$ involved the choice of some $n$-dimensional cube of group extensions. In particular, any $n$-fold extension with ``terminal vertex'' $A$ and with all other groups involved free groups, constitutes a valid choice. Note that it makes sense to call such $n$-fold extension an \emph{$n$-fold (projective) presentation} of $A$.

Let us, for one moment, return our attention to Hopf's original formula \ref{Hopfformule} and, in particular, to its denominator $[K[p],P]$. This normal subgroup of $P$ is universal in turning the projective presentation $p\colon P\to A$ into a central extension: more precisely, it is the smallest normal subgroup $N$ of $P$ with the property that the homomorphism $\overline{p}\colon P/N\to A$ induced by $p$ is a central extension. (By \emph{central extension} we mean a surjective homomorphism $f\colon A\to B$ with the property that the kernel $K[f]$ lies entirely in the center of $A$.) It was Janelidze who realized \cite{Janelidze:Hopf-talk} that the denominators of the higher Hopf formulae satisfy a similar property: they turn the higher presentations considered universally into \emph{higher central extensions}, by nulling them out of the initial vertex of the higher presentation (see also \cite{Janelidze:Hopf}). For instance, the double presentation \ref{doublepresentation} induces a double central extension
\[
\vcenter{\xymatrix{
F_{2}/I \ar[r]^-{\overline{f_2}} \ar[d]_{\overline{f'_2}} & F_{1} \ar[d]^{f_1} \\
F_{0} \ar[r]_-{f_0} & A,}}
\]
where $I$ denotes the normal subgroup $[K[f_2],K[f'_2]] \cdot [K[f_2]\cap K[f'_2], F_2]$ of $F_2$ and the homomorphisms $\overline{f_2}$ and $\overline{f'_2}$ are the natural ones induced by $f_2$ and $f'_2$, respectively.

This insight inspired Everaert, Gran and Van der Linden to write \cite{EGV}. In this article, the authors introduced notions of \emph{higher extensions} and \emph{higher central extensions}, again derived from Categorical Galois Theory, that make sense in a large variety of categories, not just in the category of groups. Making use of these notions, they obtained Hopf formulae for the cotriple homology of Barr and Beck \cite{Barr-Beck,EverVdL2}. Since Barr and Beck's theory incorporates the classical group homology,  this generalized Brown and Ellis's result to other categories than the category of groups alone. In fact, the generalization was much wider since, for instance in the case of groups, Hopf formulae were obtained not only for the integral homology, but also for the homology with coefficients in the cyclic group $\Z_n$, just to mention one example. 

One thing that is particularly interesting about the results from \cite{EGV} is that they provide one with a new approach to non-abelian (categorical) homology, based on Categorical Galois Theory or, more precisely, on the theory of higher extensions and higher central extensions, rather than on that of cotriples and simplicial objects, like in \cite{Barr-Beck}. It is therefore important that the notions of higher extension and of higher central extension be well understood. The aim of the present article is to study their fundamental properties and thus to provide a powerful base for further research in homology theory via the Brown-Ellis-Hopf formulae. 

Instead of defining higher extensions explicitly, as in \cite{EGV}, we shall here take an axiomatic approach. This will allow us to treat $n$-fold extensions (for any $n\geq 1$) as if they were ``simple'' ($1$-fold) extensions, but also, at the same time, as if they were $(n-1)$-fold extensions. This greatly simplifies the study of higher extensions and higher central extensions. Also, it will allow us to obtain a simple and direct proof of the invariance of the higher Hopf formulae with respect to the considered higher presentation, without any reference to simplicial objects or any condition on the existence of a suitable cotriple, like in \cite{EGV}.

Note that higher extensions do not explicitly appear in \cite{Brown-Ellis}. In fact, Brown and Ellis demanded their ``higher presentations'' to satisfy some weaker property instead. However, as is explained in an erratum which is available from Ellis's homepage, this was an error, and the ``higher presentations'' considered in \cite{Brown-Ellis} needed to be higher extensions after all. Donadze, Inassaridze and Porter, who gave in \cite{Donadze-Inassaridze-Porter} a new proof for Brown and Ellis's result, 
\emph{did} assume a condition on their ``$n$-fold presentations'' (or \emph{free exact $n$-fold presentations} in their terminology) which can be shown to be equivalent to demanding them to be $n$-fold extensions. 
The first appearance of both higher extensions and higher central extensions, however, was in \cite{Janelidze:Hopf-talk} although \emph{double} extensions and \emph{double central} extensions had been studied already in \cite{Janelidze:Double}. The concept of double extension has been of importance also in \cite{Carboni-Kelly-Pedicchio} and in \cite{Bourn2003} (in the latter article the term \emph{regular pushout} was used), where it was considered in a more general context than just the variety of groups, namely in \emph{regular Mal'tsev categories}. Double \emph{central} extensions in Mal'tsev \emph{varieties} have been studied in \cite{Gran-Rossi}. Finally, for $n\geq 3$, $n$-fold extensions and $n$-fold central extensions have been considered outside the scope of the variety of groups for the first time in \cite{EGV}.

The context in which the results of \cite{EGV} hold true, and in which those of the present article will be developed, is that of \emph{semi-abelian categories}. These were introduced in  \cite{Janelidze-Marki-Tholen} in order to capture, among other things, the homological properties of the categories of groups, rings, Lie algebra's, (pre-) crossed modules and similar non-abelian structures. A category $\Ac$ is called \emph{semi-abelian} if it is finitely complete and cocomplete, pointed, Barr-exact and Bourn-protomodular.

$\Ac$ being \emph{pointed} means that $0=1$, i.e.\ the initial and terminal objects coincide. This allows for a natural definition of kernels and cokernels: the \emph{kernel} of a morphism $f\colon A\to B$ is the arrow $\ker f\colon K[f]\to A$ obtained by pulling back  along $f$ the unique morphism $0\to B$; dually the cokernel of $f$ is the arrow $\coker f\colon B\to \Cok[f]$ obtained by pushing out along $f$ the unique morphism $A\to 0$. For the sake of convenience, we shall refer also to the \emph{objects} $K[f]$ and $\Cok[f]$ as to the kernel and cokernel of $f$, respectively.    

In order to explain what Barr-exactness means, let us recall that a \emph{regular epimorphism} is an arrow that is the coequalizer of some pair of arrows. Furthermore, an internal equivalence relation $(R,\pi_1,\pi_2)$ on an object $A$ in a category $\Ac$ (where $\pi_1,\pi_2\colon R\to A$ denote the projections) is called \emph{effective} if $(\pi_1,\pi_2)$ is the kernel pair of some arrow. Now, $\Ac$ being \emph{Barr-exact} \cite{Barr} means that the regular epimorphisms are pullback-stable (i.e.\ $\Ac$ is regular \cite{Barr}) and that every internal equivalence relation in $\Ac$ is effective. 

Finally, a pointed and regular category is called \emph{Bourn-protomodular} \cite{Bourn1991,Bourn2001} if the ``regular'' Short Five Lemma holds: for any commutative diagram
\[
\xymatrix{
K[f_1] \ar[r]^-{\ker f_1} \ar[d]_u & A_1 \ar[r]^{f_1} \ar[d]_v & B_1 \ar[d]^w \\
K[f_0] \ar[r]_-{\ker f_0} & A_0 \ar[r]_{f_0} & B_0}
\]  
 with $f_0$ and $f_1$ regular epimorphisms, $v$ is an isomorphism as soon as $u$ and $w$ are isomorphic.

Among the implications of the above axioms let us mention here that any regular epimorphism is the cokernel of its kernel \cite{Bourn2001}. Thus, both rows in the diagram above are short exact sequences, which means that the sequences 
\[
\xymatrix{
0\ar[r] & K[f_i] \ar[r]^-{\ker f_i} & A_i \ar[r]^{f_i} & B_i \ar[r] & 0}
\]
are exact (for $i=0,1$)  (a sequence $\xymatrix{\dots \ar[r]^{d_j} & A_j \ar[r]^{d_{j-1}} & \dots}$ is \emph{exact} if $\im d_j = \ker d_{j-1}$ for any $j$). Hence, the regular Short Five Lemma coincides with the ``classical'' Short Five Lemma. 

Let us conclude this introduction by mentioning one final property of semi-abelian categories. Recall that a monomorphism is called \emph{normal} if it is the kernel of some arrow. In a semi-abelian category, normal monomorphisms are preserved under regular images: if $k\colon K\to A$ is a normal monomorphism and $f\colon A\to B$ a regular epimorphism, then the image  $\im (f\circ k)$ of $f\circ k$ is a normal monomorphism as well \cite{Janelidze-Marki-Tholen}.

Unless stated otherwise, throughout this article, $\Ac$ will always denote a semi-abelian category.
\section{Higher extensions}

Let $\Ac$ be a semi-abelian category, with class of objects $\Ob\Ac$. If $\Ec$ is a class of morphisms in $\Ac$, then we shall write $\Ec^-$ for the class of objects $A\in\Ob\Ac$ defined as follows: $A\in\Ec^-$ if and only if there exists in $\Ec$ at least one arrow $f:A\to B$ or one arrow $g:C\to A$. The full subcategory of $\Ac$ determined by $\Ec^{-}$ will be denoted by $\Ac_{\Ec}$
\begin{definition}\label{extension}
Let $\Ec$ be a class of regular epimorphisms in $\Ac$ with $0\in\Ec^-$. We call $\Ec$ a \emph{class of extensions} when it satisfies the following list of properties. An arrow $f\in\Ec$ will be called an $\Ec$-extension, or simply an extension.
\begin{enumerate}
\item\label{axiom1}
$\Ec$ contains all split epimorphisms $f:A\to B$ with $A$ and $B$ in $\Ec^-$;
\item \label{composition}
For each composable pair of arrows $f:A\to B$ and $g:B\to C$, one has
\begin{itemize}
\item
If $f$ and $g$ are in $\Ec$, then so is $g\circ f$;
\item
If $g\circ f$ is in $\Ec$ and $B$ in $\Ec^-$, then $g$ is in $\Ec$;
\end{itemize}
\item\label{pullbackstable}
For every $f:A\to B$ in $\Ec$ and $g:C\to B$ in $\Ac$ with $C\in\Ec^-$, the pullback $g^*f$ is again in $\Ec$: 
\[
\xymatrix{
C\times_B A  \ar@{}[rd]|<<{\pullback} \ar@{.>}[r] \ar@{.>}[d]_{g^*f} & A\ar[d]^f\\
C \ar[r]_g    & B;
}
\]
\item\label{cokernel0}
If the following diagram is a short exact sequence in $\Ac$ with $A\in\Ec^-$,
\[
\xymatrix{
0\ar[r] & K\ar[r] & A \ar[r]^f & B\ar[r] & 0}
\]
then $K\in \Ec^-$ implies $f\in\Ec$;
\item\label{shortfiveregular}
For every commutative diagram with short exact rows in $\Ac$, such that the right hand vertical arrow is an identity, 
\[
\xymatrix{
0\ar[r] & K \ar[r] \ar[d]_k & A \ar[r]^f \ar[d]_a & B \ar@{=}[d] \ar[r] & 0\\
0 \ar[r] & L \ar[r] & C \ar[r]_g & B\ar[r] & 0,}
\]
one has: if both $f$ and $k$ are in $\Ec$ and $C\in\Ec^-$, then also $a\in\Ec$.
\end{enumerate}
\end{definition}

\begin{remark}\label{kernelextension} 
Since $0\in\Ec^-$,  (\ref{pullbackstable}) implies that the kernel $K[f]$ of an extension $f:A\to B$ is always in $\Ec^-$. Therefore, the implication in \ref{extension} (\ref{cokernel0}) is an equivalence.
\end{remark}
\begin{example}\label{regularexample}
The class $\Reg\Ac$ of all regular epimorphisms in $\Ac$ is a class of extensions: (\ref{composition}) and (\ref{pullbackstable}) are satisfied because $\Ac$ is regular, (\ref{shortfiveregular}) because $\Ac$ is semi-abelian: the image $\im (a)$ of $a$ is an isomorphism by the Short Five Lemma. Note that $(\Reg\Ac)^-=\Ob\Ac$. In fact, if we demand that $\Ec^-=\Ob\Ac$ then, by \ref{extension} (\ref{cokernel0}), $\Reg\Ac$  is the \emph{only} class of extensions $\Ec$ in $\Ac$: in a semi-abelian category every regular epimorphism is the cokernel of its kernel.
\end{example}

We now consider the category  $\Arr\Ac$ of all arrows in $\Ac$ which has as morphisms between arrows $a$ and $b$ commutative squares 
\[
\vcenter{\xymatrix{
A_1 \ar[r]^{f_1} \ar[d]_a &B_1 \ar[d]^b\\
A_0 \ar[r]_{f_0} & B_0}}
\]
in $\Ac$. Such a morphism will be denoted by $(f_1,f_0)\colon a\to b$.

\begin{definition}\label{doubleextension}Let $\Ec$ a be class of extensions in $\Ac$. We call \emph{double $\Ec$-extension} in $\Ac$ any arrow $(f_1,f_0)\colon a\to b$ in $\Arr\Ac$ such that all arrows in the following commutative diagram are $\Ec$-extensions, where $r$ is defined as the unique factorization to the pullback $P=A_0\times_{B_0} B_{1}$:
\[
\vcenter{\xymatrix{ 
A_{1}  \ar@/^/@{->}[drr]^{f_{1}} \ar@/_/@{->}[drd]_{a} \ar[rd]^r \\
& P\ar@{}[rd]|<{\pullback} \ar[r] \ar[d] & B_{1} \ar@{->}[d]^{b} \\
& A_0 \ar@{->}[r]_{f_0} & B_0.}}
\]
We denote the class of double $\Ec$-extensions by $\Ec^1$.
\end{definition}
Usually, we will speak of \emph{double extension}, forgetting the reference to $\Ec$. 

\begin{remark}\label{extensionsymmetry}
Of course, an arrow $(f_1,f_0)\colon a \to b$ 
in $\Arr\Ac$ is a double extension if and only if $(a,b)\colon f_1\to f_0$ is one. For \emph{higher extensions}, a similar property will be shown to hold (see Theorem \ref{higherextensionsymmetry}).
\end{remark}

\begin{remark}
Note that, because an extension is a regular epimorphism, a double extension $f$ in $\Ac$ induces a square in $\Ac$ which is a regular pushout in the sense of Bourn \cite{Bourn2003,Carboni-Kelly-Pedicchio}: every arrow in the square is a regular epimorphism and so is the factorization to the pullback. In particular, $f$ is a regular epimorphism in $\Arr\Ac$. 
\end{remark}

$\Ec^1$ is always a class of extensions in $\Arr\Ac$, as we shall now prove. \begin{lemma}\label{rotlemmalemma}
Consider in $\Ac$ a commutative diagram with short exact rows, such that $f_1$, $f_0$, $a$ and $b$ are extensions:
\[
\xymatrix{
0\ar[r] & K_1 \ar[r] \ar[d]_k & A_1 \ar[r]^{f_1} \ar[d]_a & B_1 \ar[d]^b\ar[r] & 0 \\
0 \ar[r] &K_0 \ar[r] & A_0  \ar[r]_{f_0} & B_0\ar[r] & 0.}
\]
The right hand square is a double extension if and only if $k$ is an extension.  
\end{lemma}
\begin{proof}
We can decompose the diagram as follows.
\[
\xymatrix{
0\ar[r] &K_1 \ar@{}[rd]|<<{\pullback} \ar[r] \ar[d]_k & A_1 \ar[r]^{f_1} \ar[d]_r & B_1 \ar@{=}[d] \ar[r] & 0\\
0\ar[r] &K_0 \ar[r] \ar@{=}[d] & P \ar@{}[rd]|<<{\pullback} \ar[r] \ar[d] & B_1 \ar[d]^b\ar[r] & 0\\
0\ar[r] &K_0 \ar[r] & A_0 \ar[r]_{f_0} & B_0\ar[r] &0}
\]
If the factorization $r$ to the pullback $P=A_0\times_{B_0} B_1$ is an extension, then so is $k$, by \ref{extension} (\ref{pullbackstable}) and Remark \ref{kernelextension}. If $k$ is an extension, then $r$ is one as well, by \ref{extension} (\ref{shortfiveregular}) and (\ref{pullbackstable}). 
\end{proof}

\begin{proposition}\label{doubleisextension}
$\Ec^1$ is a class of extensions in $\Arr\Ac$.
\end{proposition}
\begin{proof}
Note that $(\Ec^1)^-=\Ec$. Of course, $1_0\in\Ec$. We must prove that $\Ec^1$ satisfies the properties listed in Definition \ref{extension}. Let us refer to them as $(1)^1$, \dots, $(5)^1$, and to $(1)$, \dots, $(5)$ for the corresponding properties of $\Ec$. Keeping in mind Remark \ref{extensionsymmetry}, it is easily verified that $(\ref{axiom1})^1$ follows from $(\ref{axiom1})$ and Lemma \ref{rotlemmalemma}; $(\ref{composition})^1$ from $(\ref{composition})$ and Lemma \ref{rotlemmalemma}; $(\ref{pullbackstable})^1$ follows from $(\ref{pullbackstable})$ and $(\ref{composition})$; $(\ref{shortfiveregular})^1$ from $(\ref{pullbackstable})$ and $(\ref{shortfiveregular})$; $(\ref{cokernel0})^1$ from $(\ref{composition})$, $(\ref{cokernel0})$ and Lemma \ref{rotlemmalemma}.
\end{proof}

Proposition \ref{doubleisextension} allows us now to define, inductively, \emph{$n$-fold extensions}, for every $n\geq 1$. For this, it is necessary that we make precise what we shall mean by \emph{$n$-fold arrow}, and that we fix some notations.

Let us consider the natural numbers by their standard (von Neumann) construction and put $0=\varnothing$ and $n=\{0,\dots, n-1\}$ for $n\geq 1$. We write $\Pc(n)$ for the powerset of $n$. Recall that $\Pc(n)$ is a category which has as an arrow $S\to T$ the inclusion $S\subseteq T$ for subsets $S, T\subseteq n$.

\begin{definition}
Let $n\geq 0$. We shall call \emph{$n$-fold arrow} in $\Ac$ any contravariant functor 
\[
A\colon \Pc (n)^{\op}\to \Ac.
\]
A morphism between $n$-fold arrows $A$ and $B$ in $\Ac$ is a natural transformation $f\colon A\to B$. We will write $\Arrn\Ac$ for the corresponding category.
\end{definition}

\begin{remark}
Let $n\geq 0$. Just as any category of presheaves in $\Ac$, $\Arrn\Ac$ is a semi-abelian category. Moreover, limits (colimits) in $\Arrn\Ac$ are pointwise limits (colimits) in $\Ac$.
\end{remark}

If $A$ is an $n$-fold arrow and $S$ and $T$ are subsets of $n$ such that $S\subseteq T$, then we will write $A_S$ for the image $A(S)$ of $S$ by the functor $A$ and $a_S^T\colon A_T\to A_S$ for the image $A(S\subseteq T)$ of $S\subseteq T$. If $f\colon A\to B$ is a morphism between $n$-fold arrows, we will write $f_S\colon A_S\to B_S$ for the $S$-component of the natural transformation $f$. Furthermore, it will often be convenient to write $(A_S)_{S\subseteq n}$ instead of $A$ and $(f_S)_{S\subseteq n}$ instead of $f$, or simply $(A_S)_{S}$ and $(f_S)_S$. Moreover, in order to simplify our notations, we will write $a_i$ instead of $a^n_{n\backslash \{i\}}$, for $0\leq i \leq n-1$.

Of course, $\Arr^0\Ac\cong \Ac$. Suppose now that $n\geq 1$. An $n$-fold arrow $A$ in $\Ac$ naturally induces an arrow
\[
\Bigl(a_S^{S\cup \{n-1\}}\Bigr)_{S\subseteq n-1}\colon (A_{S\cup \{n-1\}})_{S\subseteq n-1}\to (A_{S})_{S\subseteq n-1}
\]
in $\Arr^{n-1}\Ac$. This yields a functor $\delta_{n-1} \colon\Arrn\Ac\to\Arr(\Arr^{n-1}\Ac)$. Obviously, we have
\begin{lemma}
$\Arrn\Ac\cong\Arr(\Arr^{n-1}\Ac)$. An isomorphism is given by the functor $\delta_{n-1}$.
\end{lemma}
We are now in a position to define higher extensions.

\begin{definition}
Let $\Ec$ be a class of extensions in $\Ac$. We put $\Ec^0=\Ec$ and let, inductively for any $n\geq 1$, $\Ec^n$ be the class of arrows $f\colon A\to B$ in $\Arrn\Ac$ such that the induced square in $\Arr^{n-1}\Ac$
\[
\xymatrix{
(A_{S\cup \{n-1\}})_{S}  \ar[d]_{(a_S^{S\cup \{n-1\}})_{S}}   \ar[rr]^{(f_{S\cup \{n-1\}})_{S}}& & (B_{S\cup \{n-1\}})_{S}  \ar[d]^{(b_S^{S\cup \{n-1\}})_{S}}    \\
(A_{S})_{S}  \ar[rr]_{(f_{S})_{S}}&& (B_{S})_{S}
}
\]
is a double $\Ec^{n-1}$-extension: all arrows in the above square are $\Ec^{n-1}$-extensions and so is the universal arrow to the induced pullback. Suppose $n\geq 1$. We call \emph{$n$-fold $\Ec$-extension} or  simply \emph{$n$-fold extension}, any $n$-fold arrow $A$ for which $\delta_{n-1}(A)$ is in $\Ec^{n-1}$. \end{definition}

For the sake of convenience, we shall sometimes call \emph{$0$-fold $\Ec$-extension} or simply \emph{$0$-fold extension} any object $A\in\Ec^-$.

By Proposition \ref{doubleisextension} we have
\begin{corollary}
If $\Ec$ is a class of extensions in $\Ac$, $\Ec^n$ is a class of extensions in $\Arrn\Ac$ for each $n\geq 0$.
\end{corollary}

For $n\geq 0$, denote by $\Extn\Ac$ the full subcategory of $\Arrn\Ac$ of all $n$-fold extensions. In particular, $\Ext^0\Ac=\Ac_{\Ec}$. More generally, we have that $\Extn\Ac=(\Arrn\Ac)_{\Ec^n}$, for any $n\geq 0$.

\section{Higher central extensions}\label{ce}
The categorical theory of central extensions has been developed in \cite{Janelidze-Kelly} in the context of exact Mal'tsev categories \cite{Carboni-Kelly-Pedicchio}: these are the Barr-exact categories in which every internal reflexive relation is an effective equivalence relation. In particular, the theory applies to semi-abelian categories \cite{Borceux-Bourn,Bourn1996}. It is, on the one hand,  an application of Categorical Galois Theory \cite{Janelidze:Pure}, and on the other hand a generalization of the Fr\"ohlich school's work  on central extensions in varieties of $\Omega$-groups \cite{Froehlich,Furtado-Coelho,Lue}. Where the Fr\"ohlich school's notion of central extension depended upon the choice of a subvariety of the variety of $\Omega$-groups considered, Janelidze and Kelly's notion depends on the choice of a \emph{Birkhoff subcategory} of the considered exact Mal'tsev category. In order to be able to consider \emph{higher} central extensions, we shall now generalize the notion of central extension, making it dependent on a class of extensions $\Ec$. For this, we must generalize in a similar fashion the notion of Birkhoff subcategory.  

Recall from \cite{Janelidze-Kelly} that a Birkhoff subcategory of an exact Mal'tsev category $\Ac$ is a full and replete reflective subcategory $\Bc$ of $\Ac$ that is closed in $\Ac$ under subobjects and regular quotients. 

\begin{example}
In the case where $\Ac$ is a variety of (finitary, one-sorted) universal algebras, the notion of Birkhoff subcategory coincides with that of a subvariety. 
\end{example}

\begin{remark}
As explained in \cite{Janelidze-Kelly}, the closedness of $\Bc$ under subobjects and regular quotients in $\Ac$ may equivalently be stated as follows. Let us assume that $\Bc$ is a full and replete reflective subcategory of  a regular category $\Ac$. We denote by $I$ the reflector $\Ac\to \Bc$ and by $\eta$ the unit of the adjunction
\[
\xymatrix{\Ac \ar@<1 ex>[r]^-{I} \ar@{}[r]|-{\perp} & \ \Bc. \ar@{_(->}@<1 ex>[l] }
\] 
Then $\Bc$ is closed in $\Ac$ under subobjects and regular quotients if and only if for every regular epimorphism $f:A\to B$ in $\Ac$ the square \ref{Birkhoffsquare}
below is a pushout of regular epimorphisms which, in an exact Mal'tsev category, is equivalent to it being a regular pushout (see \cite{Carboni-Kelly-Pedicchio}). 
\end{remark}

The following definition is now natural. As before, we assume that $\Ac$ is a semi-abelian category.
\begin{definition}\label{strongly}
Let $\Ec$ be a class of extensions in $\Ac$, and $\Bc$ a full and replete reflective subcategory of $\Ac_{\Ec}$. Denote by $\eta$ the unit of the corresponding adjunction and by $I\colon \Ac_{\Ec}\to\Bc$ the reflector. We call $\Bc$ a \emph{strongly $\Ec$-Birkhoff} subcategory of $\Ac$ if for every ($\Ec$-) extension $f:A\to B$ the induced square
\begin{equation}\label{Birkhoffsquare}
\vcenter{\xymatrix{
A \ar[r]^-{\eta_A} \ar[d]_f & IA \ar[d]^{If}\\
B \ar[r]_-{\eta_B} & IB}}
\end{equation}
is a double ($\Ec$-) extension.
\end{definition}
\begin{remark}
Definition \ref{strongly} first appeared in \cite{EGV}, in the cases where $\Ac$ is $\Arrn\Xc$, for $\Xc$ a semi-abelian category, $n\geq 0$ and and $\Ec$ the class of $(n+1)$-fold $\Reg\Xc$-extensions (see Example \ref{regularexample}).
\end{remark}

\begin{example}\label{regularbirkhoffexample}
In the situation where $\Ec=\Reg\Ac$, the notions of Birkhoff subcategory and of strongly $\Ec$-Birkhoff subcategory of $\Ac$ coincide. 

Let us consider the particular example where $\Ac=\Gp$ is the variety of groups and $\Ec=\Reg\Gp$. The subvariety $\Ab$ of abelian groups is a (strongly $\Reg\Gp$-) Birkhoff subcategory of $\Gp$, as is any subvariety. For any group $A$, the $A$-component of the unit of the adjunction 
\[
\xymatrix{\Gp \ar@<1 ex>[r]^-{\ab} \ar@{}[r]|-{\perp} & \  \Ab \ar@{_(->}@<1 ex>[l] }
\] 
is given by the canonical quotient $\eta_A \colon A \to \ab (A)={A}/{[A,A]}$.
\end{example}
From now on, $\Bc$ will always denote an $\Ec$-Birkhoff subcategory of $\Ac$, where $\Ac$ is semi-abelian and $\Ec$ a class of extensions in $\Ac$.

For an object $B$ of $\Ac$, let us denote $\Ext(B)$ the full subcategory of the comma category $\Ac\downarrow
B$ determined by the arrows $A\to B$ in $\Ec$.  If $g \colon C\to B$ is an arrow with $C\in \Ec^-$, we will write $g^* \colon \Ext (B)\to \Ext (C)$ for the functor that sends an extension $f\colon A\to B$ to its pullback $g^{*}f\colon C\times_B A\to C$ along $g$. By \ref{extension} (\ref{pullbackstable}), $g^*$ is well-defined. Denote the kernel pair of $f$ by $(\pi_1,\pi_2)\colon R[f]\to A$. The following are natural generalizations of Janelidze and Kelly's definitions \cite{Janelidze-Kelly}, which represent the case $\Ec=\Reg\Ac$.

\begin{definition}\label{cnt}
Let $\Ec$ be a class of extensions in $\Ac$, and $\Bc$ a strongly $\Ec$-Birkhoff subcategory of $\Ac$. Consider an extension $f\colon {A\to B}$. $f$ is
\begin{enumerate}
\item \emph{a trivial extension} (with respect to $\Bc$), when the square \ref{Birkhoffsquare} is a pullback;
\item a \emph{normal extension}, when the first projection $\pi_1 \colon R[f] \to B$ (or, equivalently, the second projection $\pi_2$) is a trivial extension;
\item a \emph{central extension}, when there exists a $g\colon{C\to B}$ in $\Ext(B)$ such that $g^{*}f\colon C\times_B A\to C$ is trivial.
\end{enumerate}
\end{definition}

\begin{remark}
It is clear that every normal extension is central. Moreover, since the pullback in $\Ac$ of a trivial extension $f:A\to B$ along an arrow in $\Ac_{\Ec}$ is again trivial, every trivial extension is normal.
\end{remark}
In our semi-abelian context, we moreover have the following. The proof from \cite{EGV} of the case $\Ec=\Reg\Ac$ remains valid.
\begin{proposition}\label{splitcentral}
Every central extension is normal. 
\end{proposition}

\begin{example}\label{regularcentralextensionsexample}
Consider again the situation of Example \ref{regularbirkhoffexample}. One easily sees that the trivial extensions are exactly the surjective homomorphisms of groups $f \colon A \to B$ with the property that the restriction $f \colon [A,A] \to [B,B]$ of $f$ to the derived subgroups is an isomorphism. It was shown in~\cite{Janelidze:Pure} that an extension $f \colon A \to B$ that is central with respect to~$\Ab$ is the same thing as a central extension in the classical sense: its kernel $K[f]$ is contained in the centre $Z(A)$ of $A$.  
\end{example}

The full subcategory of $\Ext\Ac$ of all extensions that are central with respect to $\Bc$ will be denoted by $\CExt_{\Bc}\Ac$. We have the following important property. Again, the proof from \cite{EGV} of the case $\Ec=\Reg\Ac$ remains valid.

\begin{proposition}
$\CExt_{\Bc}\Ac$ is a strongly $\Ec^1$-Birkhoff subcategory of $\Arr\Ac$.
\end{proposition}
We denote the reflector $\Ext\Ac\to \CExt_{\Bc}\Ac$ by $I_1$. In order to recall its construction, we introduce some notations.

If $I\colon \Ac_{\Ec}\to\Bc$ is the reflector associated with a strongly $\Ec$-Birkhoff subcategory $\Bc$ of $\Ac$ with unit $\eta$, we shall write $[A]_{\Bc}$ or simply $[A]$ for the kernel $K[\eta_A]$ of the $A$-component of the unit $\eta$. Moreover, we denote the arrow $\ker\eta_A\colon [A]\to A$ by $\mu_A$. This defines a functor $[\cdot]=[\cdot]_{\Bc} \colon\Ac_{\Ec}\to\Ac_{\Ec}$ together with a natural transformation $\mu\colon [\cdot]\to 1_{\Ac_{\Ec}}$. 

Now, the ``centralization'' $I_1f$ of an extension $f\colon A\to B$ is given by the extension $I_1f\colon A/[f]_{1,{\Bc}}\to B$  induced by $f$, where the normal monomorphism $[f]_{1,{\Bc}}\to A$ is obtained as the composite $\mu^1_f=\mu_A\circ [\pi_2]\circ\ker [\pi_1]$ ($(\pi_1,\pi_2)$ denotes the kernel pair of $f$):
\[
\xymatrix{
[f]_{1,{\Bc}}=K[[\pi_1]] \ar[r]^-{\ker [\pi_{1}]} \ar[d] & [R[f]] \ar[d]_{\mu_{R[f]}} \ar@<0.5 ex>[r]^-{[\pi_1]} \ar@<-0.5 ex>[r]_-{[\pi_2]} & [A] \ar[d]^{\mu_A}\\
K[\pi_1]  \ar[r]_-{\ker \pi_{1}} & R[f]  \ar@<0.5 ex>[r]^-{\pi_1} \ar@<-0.5 ex>[r]_-{\pi_2} & A,}
\]
\begin{remark}
$\mu^1_f$ is a monomorphism in $\Ac$ because so are both $\mu_A$ and $[\pi_2]\circ\ker ([\pi_1])$: $\mu_A$ by assumption, and $[\pi_2]\circ\ker ([\pi_1])$ because it is the normalization of the effective equivalence relation $([R[f]], [\pi_1],[\pi_2])$ (Note that, since $R[f]$ is clearly a reflexive relation, it is indeed an effective equivalence relation because $\Ac$ is an exact Mal'tsev category). Furthermore, since $\mu_A$ is a monomorphism, the left hand square is a pullback, hence $\mu_{R[f]}\circ \ker [\pi_1]$ is a \emph{normal} monomorphism as an intersection of normal monomorphisms. It follows that $\mu^1_f=\pi_2\circ\mu_{R[f]}\circ \ker [\pi_1]$, its regular image along $\pi_2$, is normal in $A$.  
\end{remark}

\begin{example}\label{centralization}
Consider again the situation of Examples \ref{regularbirkhoffexample} and \ref{regularcentralextensionsexample}. In this case $I_1=\ab_1$ is the reflection $\Ext\Gp \to \CExt_{\Ab} \Gp$ sending an extension $f \colon A \to B$ of groups to its centralization $\ab_1 f \colon {A}/{[K[f],A]}$ $ \to B$, where $[K[f],A]$ denotes the commutator of subgroups $K[f]$ and $A$ of $A$. 
\end{example}

The fact that $\CExt_{\Bc}\Ac$ is a strongly $\Ec^1$-Birkhoff subcategory of $\Arr\Ac$ allows us to define \emph{$2$-fold central extensions} as those double extensions that are central with respect to $\CExt_{\Bc}\Ac$ and then, inductively, to define $n$-fold central extensions, for all $n\geq 1$.

More precisely:
\begin{definition}
Put $\Bc_0=\Bc$. Inductively, for $n\geq 1$, we call an $n$-fold extension $A$ an \emph{$n$-fold central extension} if $\delta_{n-1}(A)$ is central with respect to $\Bc_{n-1}$. We write $\CExt^n_{\Bc}\Ac$ for the full subcategory of $\Ext^n\Ac$ of all $n$-fold central extensions and put $\Bc_n=\CExt^n_{\Bc}\Ac$.
\end{definition}

We shall usually say \emph{double central extension} instead of $2$-fold central extension. Also, we shall sometimes call  \emph{$0$-fold central extension} any object $A$ of the strongly $\Ec$-Birkhoff subcategory $\Bc$ of $\Ac$ and put $\CExt^0_{\Bc}\Ac= \Bc$. $I_0$ is understood to be the reflector $I\colon \Ac_{\Ec}\to \Bc$.

\begin{remark}
For any $n\geq1$, the reflector $\Extn\Ac\to \CExt^n_{\Bc}\Ac$, which we denote by $I_n$, is uniquely determined by the commutativity of the square 
\[
\xymatrix{
\Extn\Ac \ar[r]^-{I_n} \ar[d]_{\delta_{n-1}} & \Extn\Ac\ar[d]^{\delta_{n-1}} \\
\Ext(\Ext^{n-1}\Ac) \ar[r]_{(I_{n-1})_1} & \Ext(\Ext^{n-1}\Ac). }
\]
\end{remark}

For any $n\geq 0$, let us denote by $\iota^n$ the functor $\Ac\to \Arrn\Ac$ that sends an object $A$ of $\Ac$ to the $n$-fold arrow $\iota^n A$ defined by putting $(\iota^n A)_{n}=A$ and $(\iota^n A)_{S}=0$, for any strict subset $S\subsetneq n$. We write $\eta^n$ for the unit of the adjunction
\[
\xymatrix{\Extn\Ac \ar@<1 ex>[r]^-{I_n} \ar@{}[r]|-{\perp} & \ \CExt^n_{\Bc}\Ac. \ar@{_(->}@<1 ex>[l] }
\] 
From the construction of the $(I_{i-1})_1$ for $1\leq i \leq n$, and the above remark, it follows that the kernel $[\cdot]_{\Bc_n}$ of the unit $\eta^n$ factors over $\iota^n$. More precisely, $[A]_{\Bc_n}=\iota^n [A]_{n,\Bc}$ for every $n$-fold extensions $A$, for some functor $[\cdot]_{n,\Bc}\colon \Extn\Ac\to \Ac$. It is these functors  $[\cdot]_{n,\Bc}$ that provide the denominators of the higher Hopf formulae (see Section \ref{Hopfisindependent}). Note that we shall usually drop the reference to $\Bc$ and denote this functor by $[\cdot]_n$.

\begin{example}\label{doubleexample}
Consider, once again, the situation of Examples \ref{regularexample}--\ref{centralization}. It was shown by Janelidze~\cite{Janelidze:Double} that the double central extensions $(A_S)_{S\subseteq 2}$ with respect to $\Ab$ are precisely those double extensions with the property that $[K[a_0],K[a_1]]= 0$ and $[K[a_0] \cap K[a_1], A_{2}]= 0$. It is easily verified by using the properties of the commutator of groups that sending a double extension $(A_S)_S$ to the double extension $(A'_S)_S$ determined by putting 
\begin{itemize}
\item
$A'_{2}= A_{2}   /  ( \, [K[a_0],K[a_1]]  \, [K[a_0] \cap K[a_1], A_{2} ] \,  ) $    
\end{itemize}
and 
\begin{itemize}
\item
$A'_S=A_S$ for all $S\subsetneq 2$
\end{itemize}
defines a reflector $\Ext^2\Ac\to \CExt^2_{\Bc}\Ac$. Consequently,  this reflector is $\ab_2$ and 
\[
[A]_2 = [K[a_0],K[a_1]]  [K[a_0] \cap K[a_1], A_{2}].
\]
One can find explicit formulae for all $[\cdot]_n$ in this case, as well as in a few others, in \cite{EGV}. For instance, one obtains formally the same formulae for Lie algebra's over a fixed commutative ring, where the Lie bracket plays the role of commutator, or for precrossed modules, where the Peiffer commutator plays this role. 
\end{example}

\section{Symmetry}\label{symmetry}
In the previous sections, we defined notions of \emph{$n$-fold extension} and \emph{$n$-fold central extension} (for $n\geq 1$), making use of the functor $\delta_{n-1}$. This seems arbitrary, since there are clearly other ways of considering an $n$-fold arrow as an arrow in $\Arr^{n-1}\Ac$. However, these alternative ways do not induce any new notions, as we shall now see.

We continue to make the same assumptions on $\Ac$, $\Ec$ and $\Bc$  and use the same notations as before.

For $i\geq 0$, let us denote by $s_i$ the map $\N\to \N$ defined as follows:
\begin{itemize}
\item
$s_i(k)=k$ if $k < i$
\item
$s_i(k)=k+1$ if $k\geq i$.
\end{itemize}

We have the following simple lemma.
\begin{lemma}\label{simpleidentity}
If $i < j$, then $s_j\circ s_i = s_i\circ s_{j-1}$.
\end{lemma}

For a subset $S\subseteq \N$, we shall write $S^i$  for the image $s_i(S)$ of $S$ by the map $s_i$. Then, for $n\geq 1$ and $0\leq i  < n$, any $n$-fold arrow $A$ induces an arrow of $(n-1)$-fold arrows
\[
\delta_i A\colon(A_{S^i\cup\{i\}})_{S\subseteq n-1} \to (A_{S^i})_{S\subseteq n-1}.
\]
By definition, $A$ is an $n$-fold extension (resp.\ $n$-fold central extension) if $\delta_i A$ is in $\Ec^n$ (resp.\  in $\CExt_{\Bc_{n-1}}\Arr^{n-1}\Ac$), for $i=n-1$. When $n\geq 2$, the former in its turn means that $\delta_j\delta_{n-1}$ is a double $\Ec^{n-2}$-extension for $j=n-2$. The following theorems state that the same is true for \emph{any} $i$ and $j$. 
\begin{theorem}\label{higherextensionsymmetry}
Suppose $n\geq 2$. For any $0\leq i\leq n-1$ and $0\leq j\leq n-2$ and any $n$-fold arrow $A$ in $\Ac$,  the following properties are equivalent:
\begin{enumerate}
\item
$A$ is an $n$-fold extension;
\item
$\delta_iA\in\Ec^{n-1}$;
\item
$\delta_j\delta_i A$ is a double $\Ec^{n-2}$-extension.
\end{enumerate}
\end{theorem}
\begin{proof}
To prove the equivalence of $(1)$ and $(2)$, we will show that $\delta_iA\in\Ec^{n-1}$ if and only of $\delta_jA\in\Ec^{n-1}$ for any $n\geq 1$ and $0\leq i<j\leq n-1$ and any $n$-fold arrow $A$. We prove this property by induction on $n$. If $n=1$, there is nothing to prove. Suppose then that $n\geq 2$ and that the property holds for $n-1$.

Let $0\leq i< j\leq n-1$ and let $A$ be an $n$-fold arrow in $\Ac$ such that 
\[
\delta_iA\colon (A_{S^i\cup \{i\}})_S \to (A_{S^i})_S 
\] 
is in $\Ec^{n-1}$. Then, in particular, both $(A_{S^i\cup \{i\}})_S$ and $(A_{S^i})_S$ are $(n-1)$-extensions and by the induction hypothesis the middle and right hand vertical arrows in the following commutative diagram in $\Arr^{n-2}\Ac$ are elements of $\Ec^{n-2}$:  
\[
\xymatrix{
(K[\delta_iA]_{S^{j-1}\cup \{j-1\}})_S \ar[r] \ar[d] &   (A_{(S^{j-1}\cup \{j-1\})^i\cup \{i\}})_S \ar@{}[rd]|{\texttt{(i)}}\ar[r] \ar[d] & (A_{(S^{j-1}\cup \{j-1\})^i})_S \ar[d] \\
(K[\delta_iA]_{S^{j-1}})_S \ar[r] & (A_{(S^{j-1})^i\cup \{i\}})_S \ar[r] & (A_{(S^{j-1})^i})_S. }
\]
Furthermore, since $\delta_iA$ is a regular epimorphism in $\Arr^{n-1}\Ac$, the upper and lower right hand horizontal arrows are regular epimorpshims in $\Arr^{n-2}\Ac$, hence both rows are short exact sequences. Since $\delta_i A$ is an element of $\Ec^{n-1}$, its kernel $K[\delta_i A]$ is an $(n-1)$-extension. By the induction hypothesis, this implies that the left hand vertical arrow is in $\Ec^{n-2}$, hence, by \ref{extension} (\ref{cokernel0}), the right hand square $\texttt{(i)}$ is a double $\Ec^{n-2}$-extension. 

Let us have a look at another commutative diagram in $\Arr^{n-2}\Ac$:
\[
\xymatrix{
(K[\delta_jA]_{S^{i}\cup \{i\}})_S \ar[r] \ar[d] &  (A_{(S^{i}\cup \{i\})^j\cup \{j\}})_S \ar@{}[rd]|{\texttt{(ii)}} \ar[r] \ar[d] & (A_{(S^{i}\cup \{i\})^{j}})_S  \ar[d] \\
(K[\delta_jA]_{S^{i}})_S \ar[r] &  (A_{(S^{i})^{j}\cup \{j\}})_S   \ar[r] & (A_{(S^{i})^{j}})_S  }
\]
By Lemma \ref{simpleidentity}, the square $\texttt{(ii)}$ is identical to the square $\texttt{(i)}$ of the previous diagram, except that the horizontal arrows now point downwards. In particular, $\texttt{(ii)}$ is a double $\Ec^{n-2}$-extension. Since both rows are short exact sequences, this implies that the left hand vertical arrow $\delta_iK[\delta_j(A)]$ is in $\Ec^{n-2}$. By the induction hypothesis, this implies that $K[\delta_jA]$ is an $(n-1)$-extension, as well as the middle arrow $\delta_i (A_{S^j\cup \{j\}})_S$. Applying  \ref{extension} (\ref{cokernel0}) to the short exact sequence
\[
\xymatrix{
K[\delta_j A] \ar[r] & (A_{S^j\cup \{j\}})_S \ar[r] & (A_{S^j})_S,}
\]
we find that $\delta_jA\in\Ec^{n-1}$. 

Similarly, one proves that $\delta_jA\in\Ec^{n-1}$ implies $\delta_iA\in\Ec^{n-1}$.  

The equivalence of $(2)$ and $(3)$ for a certain $n$ follows from the equivalence of $(1)$ and $(2)$ for $n-1$, as one sees by applying properties \ref{extension} (\ref{pullbackstable}) and \ref{extension} (\ref{cokernel0}) of the classes of extensions $(\Ec^{n-2})^1$ and $\Ec^{n-1}$ in a similar fashion as above. Note that this latter equivalence holds also for $n=1$, in which case it holds trivially.
\end{proof}

\begin{theorem}\label{higheriscommutator}
Suppose $n\geq 1$. For any $0\leq i\leq n-1$ and any $n$-fold extension $A$ in $\Ac$, we have that \[
\iota^{n-1}[A]_{n,\Bc}=[\delta_iA]_{1,\Bc_{n-1}}
\]
and, consequently, that the following properties are equivalent:
\begin{enumerate}
\item
$A$ is an $n$-fold central extension (w.r.t.\ $\Bc$);
\item
$\delta_iA$ is a central extension (w.r.t.\ $\Bc_{n-1}=\CExt^{n-1}_{\Bc}\Ac$).
\end{enumerate}
\end{theorem}

\begin{proof}
It suffices to prove the first claim. For this, we show that the following property holds: for every $n\geq 1$ and $0\leq i< j\leq n-1$, $[\delta_iA]_{1,\Bc_{n-1}}=[\delta_jA]_{1,\Bc_{n-1}}$. Note that this is meaningful thanks to Theorem \ref{higherextensionsymmetry}. The proof is by induction on $n$. If $n=1$, there is nothing to prove. Let us then suppose that $n\geq 2$ and that the property holds for $n-1$. 
Let $0\leq i< j\leq n-1$. Consider, again, the squares $\texttt{(i)}$ and $\texttt{(ii)}$ from Proposition \ref{higherextensionsymmetry}. Recall that, by Lemma \ref{simpleidentity}, $\texttt{(i)}=\texttt{(ii)}$, except that the horizontal arrows from $\texttt{(i)}$ point downwards in $\texttt{(ii)}$. By taking kernel pairs both horizontally and vertically, we thus get the arrows  $\delta_iR[\delta_j(A)]$ and $ \delta_{j-1}R[\delta_iA]$ which have the same kernel pair, by commutativity of taking limits; we denote it by $R'=(R'_S)_S$: 

\[
\xymatrix{
 (R'_S)_S \ar@<0.5 ex>[r] \ar@<-0.5 ex>[r]  \ar@<0.5 ex>[d] \ar@<-0.5 ex>[d] &   (R[\delta_iA]_{S^{j-1}\cup \{j-1\}})_S   \ar@<0.5 ex>[d] \ar@<-0.5 ex>[d] \ar[rr]^-{\delta_{j-1}R[\delta_i(A)]} &&   (R[\delta A]_{S^{j-1}})_S  \ar@<0.5 ex>[d] \ar@<-0.5 ex>[d] \\
  (R[\delta_jA]_{S^i\cup \{i\}})_S \ar[d]_{\delta_iR[\delta_j(A)]}  \ar@<0.5 ex>[r] \ar@<-0.5 ex>[r] &  (A_{(S^i\cup \{i\})^j \cup \{j\}})_S \ar@{}[rrd]|{\texttt{(ii)}} \ar[rr] \ar[d] &&  (A_{(S^i\cup \{i\})^j})_S \ar[d] \\
  (R[\delta_jA]_{S^i})_S   \ar@<0.5 ex>[r] \ar@<-0.5 ex>[r]  & (A_{(S^i)^j\cup \{j\}})_S \ar[rr] & & (A_{(S^i)^j})_S.}
\]
Consider the left hand upper square in the diagram above. By forgetting the second projections of the kernel pairs and applying $[\cdot]_{n-2}$, we obtain the right hand lower square in the next diagram. Taking kernels horizontally and then vertically, we get an object $K$ and a commutative diagram with short exact rows and columns:
\[
\xymatrix{
K \ar[r] \ar[d] & [R[\delta_jA]]_{n-1} \ar[r] \ar[d] & [(A_{S^j\cup\{j\}})_S]_{n-1} \ar[d] \\
[R[\delta_iA]]_{n-1} \ar[r] \ar[d] & [(R'_S)_S]_{n-2} \ar[r] \ar[d] & [(R[\delta_iA]_{S^{j-1}\cup \{j-1\}})_S]_{n-2} \ar[d] \\
[(A_{S^i\cup\{i\}})_S]_{n-1} \ar[r] & [(R[\delta_jA]_{S^i\cup\{i\}})_S]_{n-2} \ar[r] & [(A_{(S^i\cup\{i\})^j\cup \{j\}})_S]_{n-2}. }
\]
Indeed, by the induction hypothesis, we have that
\[
\iota^{n-2}[R[\delta_jA]]_{n-1,\Bc} = [ \delta_{i}R[\delta_jA] ]_{1,\Bc_{n-2}}
\]
and
\[
\iota^{n-2}[R[\delta_iA]]_{n-1,\Bc} =  [\delta_{j-1}R[\delta_iA] ]_{1,\Bc_{n-2}}
\]
as well as
\[
\iota^{n-2} [(A_{S^i\cup\{i\}})_S]_{n-1,\Bc} = [ \delta_{j-1} (A_{S^i\cup\{i\}})_S]_{1,\Bc_{n-2}}
\]
and
\[
\iota^{n-2} [(A_{S^j\cup\{j\}})_S]_{n-1,\Bc} = [ \delta_i(A_{S^j\cup \{j\}})_S ]_{1,\Bc_{n-2}}
 \]
It follows that
\[
[\delta_iA]_{1,\Bc_{n-1}} = K = [\delta_jA]_{1,\Bc_{n-1}}  
\]
\end{proof}

\section{Higher central extensions of extensions}
We defined $n$-fold extensions and $n$-fold central extensions inductively, by considering $n$-fold arrows $A$ as arrows $\delta_iA$ between $(n-1)$-fold extensions. Another useful technique is to consider $n$-fold arrows as $(n-1)$-fold arrows in $\Arr\Ac$, which allows one to compare $n$-fold extensions with $(n-1)$-fold extensions of extensions and  $n$-fold central extensions with $(n-1)$-fold central extensions of extensions. These notions coincide, as we are now to show.

Suppose that $n\geq 1$. An $n$-fold arrow $A$ in $\Ac$ induces, for every $0\leq i\leq n-1$, the following $(n-1)$-fold arrow $\rho_iA$ in $\Arr\Ac$:
\[
\Bigl(a_{S^{i}}^{S^{i}\cup \{i\}} \colon A_{S^{i}\cup \{i\}}\to A_{S^{i}} \Bigr)_{S\subseteq n-1}
\]
This yields, for every $i$, a functor $\rho_{i}\colon\Arrn\Ac\to\Arr^{n-1}(\Arr\Ac)$, which has the following obvious property. 

\begin{lemma}\label{isoshift}
$\Arrn\Ac\cong\Arr^{n-1}(\Arr\Ac)$. An isomorphism is given by any of the functors $\rho_{i}$ ($0\leq i\leq n-1$).
\end{lemma}

Furthermore, by Lemma \ref{simpleidentity}, we have the following commutativity of $\rho_i$'s with $\delta_j$'s. We denote by $(\rho_i,\rho_i)$ the functor $\Arr\Arr^{n-1}\Ac\to \Arr\Arr^{n-2}\Arr\Ac$ that sends an arrow $A\to B$ in $\Arr^{n-1}\Ac$ to the induced arrow $\rho_iA\to \rho_iB$.

\begin{lemma}\label{shiftcommutes} Suppose $n\geq 2$ and $0\leq i,j\leq n-1$. If $i<j$, then $\delta_{j-1}\circ \rho_i= (\rho_i,\rho_i)\circ \delta_j$. If $i>j$, then $\delta_j\circ \rho_i=(\rho_{i-1},\rho_{i-1})\circ \rho_j$.
\end{lemma}
Using the above lemma's and Theorem \ref{higherextensionsymmetry}, we find:

\begin{lemma}\label{shiftextensionlemma}
Suppose $n\geq 1$. For any $0\leq i\leq n-1$, we have that $\rho_{i} (\Ec^{n}) =(\Ec^1)^{n-1}$.
\end{lemma}

Finally, using the three lemma's above and, once more, Theorem \ref{higherextensionsymmetry}, we get the following characterizations of higher extensions and higher central extensions. 
\begin{theorem}\label{shiftextension}
Suppose $n\geq 1$. For any $0\leq i\leq n-1$ and any $n$-fold arrow $A$ in $\Ac$, the following properties are equivalent:
\begin{enumerate}
\item
$A$ is an $n$-fold $\Ec$-extension;
\item
$\rho_iA$ is an  $(n-1)$-fold $\Ec^1$-extension.
\end{enumerate}
\end{theorem}

\begin{theorem}\label{1+n-1=n}
For any $n\geq1$ and and $0\leq i\leq n-1$, and any $n$-fold arrow $A$, $\iota^1[A]_{n}=[\rho_{i}A]_{n-1}$ and the square
\[
\xymatrix{
\Extn\Ac \ar[r]^-{I_n} \ar[d]_{\rho_{i}} & \Extn\Ac\ar[d]^{\rho_{i}} \\
\Ext^{n-1}(\Arr\Ac) \ar[r]_{(I_{1})_{n-1}} & \Ext^{n-1}(\Arr\Ac). }
\]
commutes. Consequently, the following properties are equivalent:
\begin{enumerate}
\item
$A$ is an $n$-fold central extension (w.r.t.\ $\Bc$);
\item
$\rho_iA$ is an $(n-1)$-fold central extension (w.r.t.\ $\Bc_1$).
\end{enumerate}
\end{theorem}

\section{The higher Hopf formulae}\label{Hopfisindependent}

We chose in this article to take an axiomatic approach to defining higher extensions rather than the explicit one from \cite{EGV}, in which was considered solely the leading example of higher $\Ec$-extensions where $\Ec$ consists of all regular epimorphisms of a particular semi-abelian category. A big advantage of this approach is that it allows one to treat $n$-fold extensions as ``simple'' extensions (albeit in a different semi-abelian category), which often makes life easier. Also, as we have shown in the previous section, it gives one the possibility to consider $n$-fold extensions as  \emph{$(n-1)$-fold} extensions (again in a different category), which allows one to use simple inductive arguments in various situations. We shall give in this section an important example of such a situation, and give a simple and direct proof of the independency of the construction of the higher Hopf formulae from the particular choice of ``higher presentation''. 

\begin{definition}
Let $\Ec$ be a class of extensions in $\Ac$. We say that an object $P$ of $\Ac$ is \emph{$\Ec$-projective} if for every arrow $p\colon P\to B$ in $\Ac$ and every extension $f\colon A\to B$ there exists at least one arrow $p'\colon P\to A$ such that $f\circ p'=p$.
\[
\xymatrix{
& P \ar@{.>}[dl]_{p'} \ar[d]^p \\
A \ar[r]_f & B}
\]
If $P$ is an $\Ec$-projective object, then any extension $p\colon P\to A$ is called an \emph{$\Ec$-projective presentation} of  $A$. Furthermore, we say that $\Ac_{\Ec}$ \emph{has enough $\Ec$-projective objects} if there exists for every object $A\in \Ec^-$ at least one $\Ec$-projective presentation $p\colon P\to A$.
\end{definition}
Usually, we shall say \emph{presentation} rather than \emph{$\Ec$-projective presentation}, assuming that $\Ec$ is understood. Similarly, we say \emph{projective} object, rather than \emph{$\Ec$-projective} object.

In this section, we shall always assume that $\Ac$ is a semi-abelian category, $\Ec$ a class of extensions in $\Ac$, $\Bc$ a strongly $\Ec$-Birkhoff subcategory of $\Ac$ and, furthermore, $\Ac_{\Ec}$ has enough $\Ec$-projective objects.

Let $p\colon P\to A$ be a presentation of an object $A$ of $\Ac_{\Ec}$. We write $\Delta p$ for the \emph{Hopf formula} 
\[
\frac{[P] \cap K[p]}{[p]_1}.
\]
We shall not prove here its invariance w.r.t.\ the presentation $p$ of $A$, since the proof from \cite{EverVdL1} of the case $\Ec=\Reg\Ac$ remains valid and does not take advantage of the axiomatic approach to extensions taken in the present article. Rather we shall focus on the \emph{higher} Hopf formulae, and show how their invariance follows from that of the ``classical'' Hopf formula above. The invariance of the latter explicitly means the following:    

\begin{proposition}\label{Hopf} 
Suppose $A\in\Ec^-$. Then $\Delta p \cong \Delta q$ for any two presentations $p$ and $q$ of $A$. \end{proposition}

\begin{example}\label{Hopfexample}
In the situation of Examples \ref{regularexample}--\ref{doubleexample}, $\Delta$ is the classical Hopf formula for the second homology of a group: for a presentation $p\colon P\to A$ of a group $A$, $\Delta p=([P,P]\cap K[p])/[K[p],P]$.
\end{example}

The higher Hopf formulae are constructed via \emph{higher presentations}, which we will now introduce. 
\begin{definition}
Suppose $A\in \Ec^-$. Let $n\geq 1$. We call an $n$-fold extension $P$ in $\Ac$ an \emph{$\Ec$-projective $n$-fold presentation} or simply an \emph{$n$-fold presentation} of the object $A$ if $P_0=A$ and, for every $0\neq S\subseteq n$, $P_S$ is $\Ec$-projective. 
\end{definition}

We shall also be considering $\Ec^n$-projective presentations $P\to A$ of an $n$-fold extension $A$. Therefore, the following characterization of $\Ec^n$-projective objects will be useful. One easily proves this by induction.
\begin{lemma}\label{higherprojectives}
The objects of $\Arrn\Ac$ that are projective with respect to the class of extensions $\Ec^n$ are precisely those $n$-fold arrows $P$ such that $P_S$ is $\Ec$-projective for every $S\subseteq n$.
\end{lemma}
It is then easily shown, also by induction, that the following properties hold. 
\begin{lemma}
Since $\Ac_{\Ec}$ has enough $\Ec$-projective objects, $\Extn\Ac$ has enough $\Ec^n$-projective objects.
\end{lemma}
\begin{lemma}\label{enoughnpresentations}
For every object $A\in \Ec^-$, there exists at least one $n$-fold presentation $P$ of $A$. 
\end{lemma}
Let us now show that Proposition \ref{Hopf} induces the invariance also of the \emph{higher} Hopf formulae. We use the notation $\Delta_nP$ for the \emph{$n$-th Hopf formula} 
\[
\frac{[P_n]\cap \bigcap_{i=0}^{n-1}K[p_i]}{[P]_n}.
\]
of an $n$-fold presentation $P$ of an object $A$ of $\Ac_{\Ec}$. 

\begin{example}
In the situation of Examples \ref{regularexample}--\ref{Hopfexample}, $\Delta_n$ is Brown and Ellis's Hopf formula for the $(n+1)$-st homology of a group: for an $n$-fold  presentation $P$ of a group $A$, \[
\Delta_n P=\frac{[P_n,P_n]\cap \bigcap_{i=0}^{n-1} K[p_i]}{\prod_{S\subseteq n} [\cap_{i\in S}K[p_i], \cap_{i\notin S}K[p_i]]}.\]
\end{example}

\begin{lemma}\label{inductivedelta}
$\iota^1\Delta_nP=\Delta_{n-1}(\rho_iP)$ for any $n\geq 2$, $0\leq i \leq n-1$, and any $n$-fold presentation of an object $A\in\Ec^-$.
\end{lemma}
\begin{proof}
On the one hand, $\iota^1 [P]_n = [\rho_iP]_{n-1}$, by Theorem \ref{1+n-1=n}. On the other hand, $[P_n]\cap K[p_i]=[p_i]_1$: indeed, since both $p_n\colon P_n\to P_{n\backslash i}$ and the identity $1_{P_{n\backslash i}}\colon P_{n\backslash i}\to P_{n\backslash i}$ are presentations of $P_{n\backslash i}$ (since $P_{n\backslash i}$ is projective), we have, by Theorem \ref{1+n-1=n}, that
\[
\frac{[P_n]\cap K[p_i]}{[p_1]_1} = \Delta p_i \cong\Delta 1_{P_{n\backslash i}} =\frac{[P_{n\backslash i}]\cap K[1_{P_{n\backslash i}}]}{[1_{P_{n\backslash i}}]}= 0.
\]
\end{proof}

\begin{theorem}\label{higherhopf} 
Let $\Ac$, $\Ec$ and $\Bc$ be as above, and $A\in\Ec^-$. Then  $\Delta_nP\cong \Delta_n Q$ for any two $n$-fold presentations $P$ and $Q$ of $A$.
\end{theorem}
\begin{proof}
The proof is by induction on $n$. Proposition \ref{Hopf} provides the case $n=1$. The induction step will follow from Lemma \ref{inductivedelta}: indeed, let us suppose that $n\geq 2$ and that the theorem is valid for all $k<n$, and for any $\Ac$, $\Ec$ and $\Bc$. Let $A$ be an object of $\Ac_{\Ec}$ and $P$ and $Q$ be $n$-fold presentations of $A$. Consider the induced $(n-1)$-fold presentations $\rho_0P$ and $\rho_1Q$ and in particular the ''presented extensions'' $(\rho_0P)_0=p_0^1$ and $(\rho_1Q)_0= q_0^{\{1\}}$. Since $\Ac_{\Ec}$ has enough projectives, we can find a double presentation of $A$ of the form
\[
\xymatrix{
R_1 \ar[r]^r \ar[d]_{r'} & P_1 \ar[d]^{p^1_0} \\
Q_{\{1\}} \ar[r]_{q_0^{\{1\}}} & A.}
\]
Moreover, since $\Ext^2\Ac$ has enough projectives as well, and by taking into account Theorem \ref{shiftextension}, we can extend this double presentation to an $n$-fold presentation $R$ of $A$, in such a way that $r_1^2=r$, $r^2_{\{1\}}=r'$, $r_0^1=p_0^1$ and $r_0^{\{1\}}=q_0^{\{1\}}$. We have that $(\rho_0P)_0=(\rho_0R)_0$ and $(\rho_1Q)_0=(\rho_1R)_0$ hence, by the induction hypothesis,  
\[
\Delta_{n-1}(\rho_0P)\cong \Delta_{n-1}(\rho_0R)
\]
and 
\[
\Delta_{n-1}(\rho_1Q)\cong \Delta_{n-1}(\rho_1R).
\]
Using Lemma \ref{inductivedelta}, we find that
$\Delta_nP \cong \Delta_nR \cong \Delta_nQ$.
\end{proof}

\begin{remark}
In \cite{EGV} the invariance of the higher Hopf formulae was obtained indirectly, as a consequence of the equivalence with Barr and Beck's cotriple homology. In order for this equivalence to hold, the existence of a ``suitable'' cotriple on the considered category $\Ac$ had to be assumed. This assumption has disappeared here, but not entirely. It has been replaced by the weaker assumption that $\Ac$ has enough projectives. Of course, if one wants to develop a homology theory based on the theory of higher central extensions, using the higher Hopf formulae as definition, it does not make sense to make any assumption on the existence of some cotriple. Therefore, a direct prove of the invariance was necessary.
\end{remark}

\begin{remark}
There is the following beautiful description of the Hopf formula which we learned from Tim Van der Linden: $\Delta p$ is the ``difference'' between the centralization and the trivialization of the projective presentation $p$. 

Indeed, the full subcategory $\TExt_{\Bc}\Ac$ of $\Ext\Ac$ of trivial extensions (w.r.t.\ $\Bc$) is reflective, the \emph{trivialization} $Tf=T_{\Bc}f$ of an extension $f\colon A\to B$ being the pullback $\eta^*_B(If)$ of $If\colon IA\to IB$ along $\eta_B\colon B\to IB$. Moreover, it is easily verified that $\iota\Delta p$ is the kernel of the factorization $I_1p\to Tp$, which is a regular epimorphism by the strongly $\Ec$-Birkhoff property of $\Bc$.

Similarly, for $1\leq i\leq n-1$, we have that $\iota^n\Delta_nP$ is the kernel of  the factorization $I_{1,\Bc_{n-1}}\del_iP\to T_{\Bc_{n-1}}\del_iP$. It follows that $T_nP=T_{\Bc_{n-1}}\del_iP$ is independent of the choice of $i$, so that the trivialization $T_n$ of $n$-fold \emph{presentations} is well defined. Hence we may write that
\[
\iota^{n-1}\Delta_nP= K[ I_nP\to T_nP],
\]
for any $n\geq 1$ and any $n$-fold presentation $P$.

Let us also remark that, in general,  $T_{\Bc_{n-1}}\del_iA$ \emph{does} depend on the choice of $i$, when $A$ is an arbitrary $n$-fold \emph{extension}: trivialization is not well defined for arbitrary $n$-fold extensions. 
\end{remark}


\providecommand{\bysame}{\leavevmode\hbox to3em{\hrulefill}\thinspace}
\providecommand{\MR}{\relax\ifhmode\unskip\space\fi MR }
\providecommand{\MRhref}[2]{%
  \href{http://www.ams.org/mathscinet-getitem?mr=#1}{#2}
}
\providecommand{\href}[2]{#2}

\noindent Tomas Everaert  \  teveraer@vub.ac.be \\
Vakgroep Wiskunde, Vrije Universiteit Brussel, Pleinlaan 2, 1050 Brussel, Belgium

\end{document}